\documentclass[11pt]{article}
\usepackage{amsmath}
\usepackage{amsthm}
\usepackage{amssymb}

\newtheorem{theo}{Theorem}[section]

\newtheorem{lemma}[theo]{Lemma}

\newcommand{\TT}{{\bf T}}
\newcommand{\FF}{{\cal F}}
\newcommand{\GG}{{\cal G}}

\newcommand{\eps}{{\varepsilon}}
\begin{document}
\date{}
\title{
Implicit representation of sparse hereditary families}
\author{Noga Alon
\thanks
{Princeton University,
Princeton, NJ 08544, USA and
Tel Aviv University, Tel Aviv 69978,
Israel.
Email: {\tt nalon@math.princeton.edu}.
Research supported in part by
NSF grant DMS-1855464 and BSF grant 2018267.}
}
\maketitle
\begin{abstract}
For a hereditary family of graphs $\FF$, let $\FF_n$ denote the set
of all members of $\FF$ on $n$ vertices. The speed of $\FF$ is the
function $f(n)=|\FF_n|$. An implicit representation of size
$\ell(n)$ for $\FF_n$ 
is a function assigning a label of $\ell(n)$ bits
to each vertex of any given graph $G \in \FF_n$, so that the adjacency
between any pair of vertices can be determined by their labels. 
Bonamy, Esperet, Groenland and Scott proved that the minimum
possible size of an implicit representation of $\FF_n$ for any
hereditary family $\FF$ with speed $2^{\Omega(n^2)}$ is 
$(1+o(1)) \log_2 |\FF_n|/n~(=\Theta(n))$.
A recent result of Hatami and Hatami shows that the situation is
very different for very sparse hereditary families. They showed
that for every
$\delta>0$ there are 
hereditary families of graphs with speed $2^{O(n \log n)}$ that do
not admit implicit representations of size smaller than
$n^{1/2-\delta}$.  In this
note we show that even a mild speed
bound ensures an implicit representation of size $O(n^c)$ for some 
$c<1$.  Specifically we prove that for every 
$\eps>0$ there is an integer $d \geq 1$ so that 
if $\FF$ is a hereditary family
with speed $f(n) \leq 2^{(1/4-\eps)n^2}$ then $\FF_n$ admits an 
implicit representation of size $O(n^{1-1/d} \log n)$. Moreover,
for every integer $d>1$ there is a hereditary family for which this is
tight up to the logarithmic factor. 
\end{abstract}
\section{Introduction}
A family of graphs $\FF$ is hereditary if it is closed under taking
induced subgraphs. Let $\FF_n$ denote the set of all members of
$\FF$ with $n$ vertices. The speed of $\FF$ is the function $f(n)=
|\FF_n|$. An implicit representation of size $\ell(n)$ of
$\FF_n$ is a function assigning a label of $\ell(n)$ bits
to
each vertex of any given graph $G \in \FF_n$, so that the adjacency
between any pair of vertices can be determined by their labels. It
is easy and well known (see \cite{KNR}) that the existence of such a
function is equivalent to the existence of a graph on
$2^{\ell(n)}$ vertices  which contains every member of $\FF_n$
as an induced subgraph (here we do not assume that the function 
assigning labels has to be efficiently computable). Such a graph is
called a universal graph for $\FF_n$. To see the equivalence observe
that given a function corresponding to an implicit representation
of size $\ell(n)$ the graph whose vertices are all possible
labels in which two are adjacent iff the corresponding labels 
determine adjacency in a graph of $\FF_n$ is a universal graph for
$\FF_n$. The converse follows by assigning to each vertex of 
a graph $G  \in \FF_n$ the number of the vertex of the universal
graph that plays its role in a copy of $G$ in this graph.

There is a vast literature dealing with universal graphs for
various families, see, e.g., \cite{Alo}, \cite{BEGS}, \cite{HH} and
the many references therein.  By the above remark, the minimum
possible size $\ell(n)$ of labels for a family $\FF_n$ has to
satisfy $[2^{\ell(n)}]^n \geq |\FF_n|$, that is, $\ell(n) \geq 
\frac{\log_2 |\FF_n|}{n}$, and it is known that this is essentially
tight in many interesting cases. In particular, this is the case 
for the family of
all graphs (see \cite{Mo}, \cite{Alo}). It is also nearly
tight for many additional examples, including all hereditary
families
satisfying  $|\FF_n| =2^{\Omega(n^2)}$. By known results
\cite{Ale}, \cite{BT}, if $|\FF_n| =2^{\Omega(n^2)}$ then
$|\FF_n| =2^{(1-1/k) n^2/2+o(n^2)}$ for some integer $k > 1$.
Bonamy, Esperet, Groenland and Scott
\cite{BEGS} proved that in all these cases there is an
implicit representation with labels of length 
$(1-1/k) n/2+o(n)$.  On the other hand, a recent result
of Hatami and Hatami \cite{HH}, settling a problem raised 
by Kannan, Naor and Rudich
\cite{KNR}, shows that there are very sparse hereditary families
for which  any implicit representation requires labels of size nearly
$\sqrt n$. Specifically it is shown in \cite{HH} that for every 
$\delta>0$ there is a hereditary family $\FF$ satisfying
$|\FF_n|=2^{O(n \log n)}$ so that the size of any implicit
representation 
for $\FF_n$ is at least $\Omega(n^{1/2-\delta})$.  It is not clear
if the exponent $1/2$ can be improved, and it is also not known
what happens for families $\FF$ with speed $f(n)$ exceeding
$2^{n \log n}$ which is $2^{o(n^2)}$.  It is known that in this
range the speed is at most $2^{n^{2-\eps}}$ for some fixed
$\eps>0$ (see \cite{ABBM}). Our contribution here is to show
that in all these cases there is an implicit representation of size
at most $O(n^{1-\eps})$.
\begin{theo}
\label{t11}
For any $\eps>0$ there is an integer $d \geq 1$ so that the
following holds.
Let $\FF$ be a hereditary family of graphs with speed
$f(n) =|\FF_n| \leq 2^{(1/4-\eps)n^2}$ (and hence 
$f(n)=2^{o(n^2)}$). Then 
there is an implicit representation of size at most
$O(n^{1-1/d} \log n)$  for $\FF_n$.  In addition, for any such
integer $d>1$ there is a hereditary family for which this is
tight up to the $\log n$ factor.
\end{theo}
Natural examples of hereditary families $\FF$ of graphs are
intersection graphs of geometric objects of prescribed type. 
In many of these cases it is possible to obtain tight bounds for
the function $f(n)=|\FF_n|$ using a theorem of Warren \cite{Wa}
from real algebraic geometry.
This theorem, as well as a related earlier one by Milnor \cite{Mi},
have been applied by Goodman and Pollack in order to estimate the
number of configuration and polytopes in $R^d$. Their results
appear in the very first volume of the journal Discrete and
Computational Geometry they founded in the mid. 80s \cite{GP}.
See also \cite{Al00}, \cite{Al0} and the brief discussion in
Section 3 here for more about this topic. 

\section{Proof}
For any two integers $k,d \geq 1$ let $U(k,d)$ denote the bipartite graph
with two vertex classes $A,B$ satisfying $|A|=d$, and $|B|=k \cdot
2^d$,
where for each subset $C \subset A$ there are exactly $k$ vertices
in $B$ whose set of neighbors in $A$ is exactly $C$. If $X,Y$ are
disjoint sets of vertices of a graph $G$, let $G[X,Y]$ denote the
bipartite graph induced by the sets $X$ and $Y$ (ignoring the edges 
inside $X$ and inside $Y$). Call a graph $U(k,d)$-free if it contains
no two disjoint sets of vertices $X,Y$ so that $G[X,Y]$ is a copy 
of $U(k,d)$. Note that the graph $U(d,d)$ contains every bipartite graph
with two classes of vertices, each of size $d$, as an induced subgraph.
Therefore, if a graph contains a copy of $U(d,d)$ then it contains at
least $2^{d^2}$ distinct labelled induced subgraphs on 
$2d$ vertices.   It thus follows that if the speed of a hereditary family
$\FF$  satisfies $f(n) \leq 2^{(1/4-\eps)n^2}$ for some fixed
$\eps>0$ then there is a finite $d=d(\eps)$ so that every graph in the
family is $U(d,d)$-free. We proceed to show that the family of 
all $U(d,d)$-free graphs admits an implicit representation of size 
at most $O(n^{1-1/d} \log n).$

A set $I$ of coordinates is shattered by a family of 
binary vectors if
the projections of these vectors on $I$ includes all
$2^{|I|}$ possible binary vectors of length $|I|$.

We need the
following lemma.
\begin{lemma}
\label{l21}
Let $\TT$ be a family of at least
$$1+(k+d-1)\cdot 2^d \cdot {t \choose d}+\sum_{i=0}^{d-1} {t \choose
i}$$ 
distinct binary vectors of  length $t$. Then there is a set $I$ of
$d$ coordinates shattered $k+d$ times, namely,
every binary function from $I$ to
$\{0,1\}$ is a projection of at least $k+d$ distinct vectors in
$\TT$ on $I$.
\end{lemma} 
\vspace{0.2cm}

\noindent
{\bf Proof:}\, As long as $\TT$ contains more than 
$\sum_{i=0}^{d-1} {t \choose i}$ vectors there is a shattered set
of $d$ coordinates, by the Sauer-Perles-Shelah Lemma 
\cite{Sa}. Removing the $2^d$ shattering vectors from $\TT$ and
repeating the argument $(k+d)$ times we get, by the pigeonhole
principle, the same $d$-set shattered $k+d$ times. \hfill $\Box$

For a binary vector $v$ let $c(v)$ denote the number of indices $i$
so that $v_i \neq v_{i+1}$. Note that these indices partition
the set of all indices into $c(v)+1$ intervals, so that $v$ is constant
on each interval.
The primal shatter function of a family of binary vectors is the
function $g(t)$ whose value is the largest number of distinct
projections of the vectors on a set of $t$ coordinates. 
The
following lemma is proved in \cite{We} (after its optimization in
\cite{Ha}), see also \cite{CW}, \cite{MWW}. The formulation 
in these references is in terms of the notion of spanning trees with 
low crossing number. The (equivalent) formulation we use here appears in
\cite{AMY}.
\begin{lemma}
\label{l22}
Let $\GG$ be a family of binary vectors of length $n$ with primal
shatter function  $g(t) \leq ct^d$ for some constant $c>0$ and 
integer $d
\geq 1$. Then there is a fixed permutation of the coordinates of the
vectors so that for each permuted vector $v$,
$c(v) \leq O(n^{1-1/d})$.
\end{lemma}
\vspace{0.2cm}

\noindent
{\bf Proof of Theorem \ref{t11}:}\,
Let $\FF$ be a hereditary family with speed $f(n) \leq
2^{(1/4-\eps)n^2}$.
By the assumption and the remark in the first paragraph of this
section there is a finite integer $d \geq 1$ so that every member
of $\FF_n$ is $U(d,d)$-free. 
For a graph $G \in \FF_n$ let
$\GG$ be the set of rows of the adjacency matrix of $G$. These are
binary vectors of length $n$. We claim that the primal shatter
function of these family of vectors satisfies $g(t) \leq 10 t^d$
for all $t>d$. Indeed, otherwise by Lemma \ref{l21} with $k=d$
there is a set $I$ of $d$-coordinates which is shattered 
$2d$ times by these vectors. This gives a set $A$ of $d$ vertices
of $G$ and another set $B'$ of $2d \cdot 2^d$ vertices so that for
every subset $C$ of $A$ there are $2d$ vertices in $B'$ whose
set of neighbors in $A$ is exactly $C$. Let $B$ be  a subset
of $B'-A$ containing exactly $d$ vertices for each such subset $C$.
This gives a copy of $U(d,d)$ contradicting the fact that $G$ contains no 
such copy. This proves the claim. Therefore by Lemma \ref{l22} there
is a numbering of the vertices so that according to this
numbering the set of all neighbors of each vertex consists of
at most $O(n^{1-1/d})$ intervals. Assign to each vertex a label
consisting of its number and the endpoints of the corresponding
intervals. This is clearly a valid implicit representation,
establishing the required upper bound.

The (near) tightness follows by using the projective norm graphs
described in \cite{ARS}. These are graphs on $n$ vertices with
$\Omega(n^{2-2/d})$ edges that contain no copy of the complete
bipartite graph $K_{d,k}$ with $k=(d-1)!+1$. Our hereditary family
$\FF$ consists of all these graphs (for all values of $n$ for which they
exist) and all their (not necessarily induced) subgraphs. This is a
hereditary family, in fact even a monotone one. It does not contain 
an induced copy of $U(k,d)$ and hence, by the argument above
which works for $U(k,d)$ just as done for $U(d,d)$,
admits an implicit representation
of size $O(n^{1-1/d} \log n)$. Here, in fact, there is a simpler
way to get the existence of such an implicit representation. By the
K\H{o}v\'ari-S\'os-Tur\'an theorem \cite{KST} every graph in
$\FF_n$
is $p=O(n^{1-1/d})$-degenerate, hence there is an ordering of the
vertices so that every vertex has at most $p$ neighbors following
it. One can thus assign to each vertex a label consisting of its
number in this ordering and the numbers of its neighbors following
it to get the required representation.
On the other hand the speed of $\FF$
satisfies $f(n) \geq 2^{\Omega(n^{2-1/d})}$ for every $n$
for which our family contains one of the projective norm graphs.
Therefore each implicit representation for $\FF_n$ requires labels
of length at least $\log |\FF_n|/n =\Omega(n^{1-1/d})$.
This completes the proof.  \hfill $\Box$

\section{Problem}
By Theorem \ref{t11} if $\FF$ is a hereditary family with speed
$f(n)=2^{o(n^2)}$ then $\FF_n$ admits an implicit representation of
size at most $O(n^{1-1/d} \log n)$ for some integer $d \geq 1$. 
It would be interesting to decide if tighter bounds hold when 
the growth rate of the speed $f(n)$ is slower. A particularly interesting
case is $f(n) \leq 2^{O(n \log n)}$, as this holds for many
interesting hereditary families including all the ones in which
every vertex is a point in a real space of bounded dimension, and the
adjacency of two vertices 
is determined by the signs of a finite  set of bounded
degree polynomials in the coordinates of the corresponding points. Such
families, which are hereditary by definition, include many
intersection graphs of
simple geometric objects of a prescribed shape. 
By a theorem of Warren from real
algebraic geometry that deals with sign patterns of
real polynomials \cite{Wa} the speed of any such family is
at most $2^{O(n \log n)}$. The argument, which is similar to 
the one given by Goodman and Pollack in \cite{GP}, found
a significant number of applications following their work.
See \cite{Al0} and the references therein for several early examples.
However, there are quite a few families of
this type for
which the existence of economic implicit representations 
is not known. Simple examples include intersection graphs of
segments or discs in the plane studied in \cite{MM}.

By the main result of \cite{HH} for any $\delta>0$ 
there are hereditary families with
speed $f(n) \leq 2^{O(n \log n)}$ so that $\FF_n$ does not admit 
an implicit representation of size smaller than
$n^{1/2-\delta}$, and the authors of \cite{HH} raise the
natural question if the constant $1/2$ can be improved. 
Is it possible that such families always admit
an implicit representation of size $O(n^{1/2} \log n)$? 
Similarly, 
if the speed is smaller than $2^{n^{1+\eps}}$ for a sufficiently
small fixed $\eps>0$, is there always an implicit representation of
size at most $O(n^{2/3} \log n)$?


\begin{thebibliography}{99}
\bibitem{Ale} V.E.~Alekseev, 
On the entropy values of hereditary
classes of graphs, Discrete Math. Appl.  3
(1993), 191--199.
\bibitem{Al00}
N. Alon,
The number of polytopes, configurations and real matroids,
Mathematika 33 (1986), 62-71.
\bibitem{Al0}
N. Alon, Tools from higher algebra, in : Handbook of
Combinatorics, R.L. Graham, M. Gr\"otschel and L. Lov\'{a}sz, eds,
North Holland (1995), Chapter 32, pp. 1749-1783.
\bibitem{Alo}
N. Alon, Asymptotically optimal induced universal graphs,
Geometric and Functional Analysis 27 (2017), 1-32.
\bibitem{ABBM}
N. Alon, J. Balogh, B. Bollob\'as and R. Morris,
The structure of almost all graphs in a hereditary
property,
J. Comb. Theory, Ser. B 101 (2011), 85--110.
\bibitem{AMY}
N. Alon, S. Moran and A. Yehudayoff,
Sign rank, VC dimension and spectral gaps,
Proc. COLT 2016, 47--80.
Also: Mat. Sbornik 208:12 (2017), 1724-1757.
\bibitem{ARS}
N. Alon, L. R\'onyai and T. Szab\'o,
Norm-graphs: variations and applications,
J. Combinatorial Theory, Ser. B 76 (1999), 280--290.
\bibitem{BEGS}
M. Bonamy, L. Esperet, C. Groenland and A. Scott,
Optimal labelling
schemes for adjacency, comparability, and reachability, 
Proc. 53rd
Annual ACM SIGACT Symposium on Theory of Computing (STOC), pages
1109--1117, 2021.
\bibitem{BT}  B.~Bollob\'as and A.~Thomason, Hereditary and
monotone properties of graphs, in: The Mathematics of Paul
Erd\H{o}s, II (R.L. Graham and J. Ne\v{s}et\v{r}il, Editors),
Alg. and Combin., Vol. 14, Springer-Verlag, New
York/Berlin (1997), 70--78.
\bibitem{CW}
B. Chazelle and E. Welzl,
Quasi-optimal range searching in
spaces of finite VC-dimension,
Discrete Comput. Geom. 4, no.  5 (1989), 467--489. 
\bibitem{GP}
J. E. Goodman and R. Pollack,
Upper bounds for configurations and polytopes in $R^d$
Discrete Comput. Geom. 1 (1986), 219--227.
\bibitem{Ha}
D. Haussler,
Sphere packing numbers for subsets of the Boolean n-cube
with
bounded Vapnik-Chervonenkis dimension,
J. Combin. Theory Ser. A 69, no. 2 (1995), 217--232.
\bibitem{HH}
H. Hatami and P. Hatami,
The implicit graph conjecture is false, arXiv:2111.13198, 2021.
\bibitem{KNR}
S. Kannan,
M. Naor and S. Rudich,
Implicit representation of graphs,
SIAM J. Discrete Math. 5 (1992), 596--603.
\bibitem{KST} 
T. K\H{o}v\'ari, V. T. S\'os and P. Tur\'an,
On a problem of K. Zarankiewicz,
Colloquium Math. 3, (1954), 50-57.
\bibitem{MM}
C. McDiarmid and T. M\"uller,
Realizations of disk and segment graphs, J. Combin. Theory Ser. B
103 (2013), no. 1, 114--143.
\bibitem{MWW}
J. Matou\v{s}ek, E. Welzl and L. Wernisch,
Discrepancy and
approximations for bounded VC-dimension,
Combinatorica 13, no. 4 (1993), 455--466. 
\bibitem{Mi}
J. Milnor,
On the Betti numbers of real varieties, Proc. AMS 15 (1964),
275--280.
\bibitem{Mo}
J. W. Moon, On minimal n-universal graphs, Proceedings of the
Glasgow Mathematical Association,
7(1) (1965), 32--33.
\bibitem{Sa} 
N.~Sauer,
On the density of families of sets,
J. Combinatorial Theory, Ser. A 13 (1972), 145--147.
\bibitem{Wa}
H. E. Warren,
Lower Bounds for approximation by nonlinear
manifolds, Trans. Amer. Math. Soc. 133 (1968), 167-178.
\bibitem{We}
E. Welzl,
Partition trees for triangle counting and other range
searching problems,
Proc. 4th Annual Symposium on
Computational Geometry,
pages 23--33, 1988.
\end{thebibliography}
\end{document}